\documentclass{amsart}%
\usepackage{amssymb}
\usepackage{amsmath}
\usepackage{amsfonts}
\usepackage{graphicx}
\usepackage{hyperref}%
\providecommand{\U}[1]{\protect\rule{.1in}{.1in}}
\newtheorem{theorem}{Theorem}
\theoremstyle{plain}

\newtheorem{lemma}{Lemma}

\numberwithin{equation}{section}
\providecommand{\U}[1]{\protect\rule{.1in}{.1in}}
\begin{document}
\title[On Concatenations of Two $ k $-Generalized Fibonacci Numbers]{On Concatenations of Two $ k $-Generalized Fibonacci Numbers}

\author{Alaa Altassan}
\address{Department of Mathematics, King Abdulaziz University, 21589, Jeddah, Saudi Arabia}
\email{aaltassan@kau.edu.sa}

\author{Murat Alan}
\address{Department of Mathematics, Yildiz Technical University, Davutpasa Campus,
Esenler, 34210 Istanbul, Turkey.}
\email{alan@yildiz.edu.tr}

\begin{abstract}
Let $ k \geq 2 $ be an integer. One of the generalization of the classical Fibonacci sequence is defined by the recurrence relation
$ F_{n}^{(k)}=F_{n-1}^{(k)} + \cdots + F_{n-k}^{(k)}$ for all $ n \geq 2$ with the initial values $ F_{i}^{(k)}=0 $ for $ i=2-k, \ldots, 0 $ and $ F_{1}^{(k)}=1.$ $  F_{n}^{(k)} $ is an order $ k $ generalization of the Fibonacci sequence and it is called $ k- $generalized Fibonacci sequence or shortly $ k- $Fibonacci sequence. In 2020, Banks and Luca, among other things, determined all Fibonacci numbers which are concatenations of two Fibonacci numbers. In this paper, we consider the analogue of this problem and determine all  $ k-$generalized Fibonacci numbers which are concatenations of two terms of the same sequence. We completely solve this problem for all $ k \geq 3. $ 

\bigskip


\bigskip


\end{abstract}
\maketitle

\section{Introduction}
Let $ k \geq 2 $ be an integer. One of the generalization of the classical Fibonacci sequence is a sequence defined by the recurrence relation
$ F_{n}^{(k)}=F_{n-1}^{(k)} + F_{n-2}^{(k)} + \cdots + F_{n-k}^{(k)}$ for all $ n \geq 2$ with the initial conditions $ F_{i}^{(k)}=0 $ for $ i=2-k, \ldots, 0 $ and $ F_{1}^{(k)}=1.$ $  F_{n}^{(k)} $ is an order $ k $ generalization of the Fibonacci sequence and it is called the $ k- $generalized Fibonacci sequence or shortly $ k- $Fibonacci sequence. Fibonacci sequence is a special case of this sequence which is corresponding the case $ k=2 .$ In 2000, Luca and Banks, among other things, determined all Fibonacci (Lucas) numbers that are concatenations of two Fibonacci (Lucas) numbers \cite{Banks}. In particular, they proved that $ F_{10}= 55$ is the greatest Fibonacci number with this property. In fact, they obtained  more general results related to finiteness of these concatenations under some mild hypotheses. In \cite{Alan}, the author found all Fibonacci (Lucas) numbers that are concatenations of two Lucas (Fibonacci) numbers whereas in \cite{Altassan} the authors take into account the mixed concatenations of these sequences. In \cite{Bravo23,Erduvan,Guneymerve}, the authors consider the concatenations of Padovan and Perrin numbers. Recently, some authors extended the investigation of concatenations of some recurrence sequences by taking into account their base$ -b $ representation \cite{Adedji1,Adedji2,Adedji3}.

In this paper, for arbitrary $ k \geq 2 $, we consider the concatenations of two $ k-$generalized Fibonacci numbers that forms another $ k-$generalized Fibonacci number. As a mathematical formulation of this problem, we consider the Diophantine equation
\begin{equation}
F_n^{(k)}  = F_m^{(k)} 10^{d} + F_l^{(k)} , \quad  k \geq 3
\label{Fnml}
\end{equation}
in non negative integers $n, m>0 $ and $ l$ where $ d $ is the number of digits of $ F_l^{(k)} .$ We state the main result of this paper as follows. 
\begin{theorem}
\label{main1}
The Diophantine equation \eqref{Fnml} has solutions only in the cases 
\begin{eqnarray*} & F_{7}^{(3)}=24= \overline{F_{3}^{(3)}F_{4}^{(3)}}  ,\\
 & F_{8}^{(3)}=44= \overline{F_{4}^{(3)}F_{4}^{(3)}}  ,\\
 & F_{16}^{(8)}=16128= \overline{F_{6}^{(8)}F_{9}^{(8)}}  
\end{eqnarray*}
\end{theorem}
It is known that for $ k=2, $ the equation in \eqref{Fnml} has solutions only for $ F_n \in \{ 13, 21, 55 \} $ \cite{Banks} so that we state Theorem \ref{main1} only for $ k \geq 3. $

There exists some effective methods for handling Diophantine equations similar to \eqref{Fnml} such as linear forms in logarithms of algebraic numbers due to Matveev \cite{Matveev} or some versions of reduction algorithm due to Dujella and Peth\H o  \cite{DP}. We use these tools together with some properties of $k-$Fibonacci numbers and continued fractions. It is worth to note that when carrying out these methods,  we used the software Maple for all computations by taking into account 1000 digit. We give the details of these methods in the next chapter.

\section{Preliminaries}

Let $\xi$ be an algebraic number and let 
\[
\sum_{j=0}^{n} u_jx^{n-j}=u_0 (x-\xi^{(1)}) \cdots (x-\xi^{(n)}) \in \mathbb{Z}
\]
be its minimal polynomial over the rational numbers where the $u_j$'s are coprime integers with $u_0>0$ and the $\xi^{(i)}$'s are conjugates of $\xi$.

The logarithmic height of $\xi$ is defined by
\[
h(\xi)=\frac{1}{n}\left(\log u_0+\sum_{i=1}^{n}\log\left(\max\{|\xi^{(i)}|,1\}\right)\right).
\]
In particular, for a rational number $ r/s,$ with $ \gcd(r, s)=1 $ and $ s>0 ,$  $h(r/s)=\log \max \{|r|,s\}$. We use the following properties in the more complicated case.
\begin{itemize}
\item[$ \bullet $] $h(\xi_1 \pm \xi_2)\leq h(\xi_1) + h(\xi_2)+\log 2$.
\item[$ \bullet $] $h(\xi_1 \xi_2^{\pm 1})\leq h(\xi_1)+h(\xi_2)$.
\item[$ \bullet $] $h(\xi^{s})=|s|h(\xi),$ $ s \in \mathbb{Z} $.
\end{itemize}

\begin{theorem}[\cite{Matveev}]
\label{Matveev}
Consider $\alpha_1, \ldots, \alpha_t$ as positive real algebraic numbers within a real algebraic number field $\mathbb{F}$ of degree $ d_\mathbb{F} ,$ and  $b_1,\ldots,b_t$ be rational integers. Define
\[
U:=\alpha_1^{b_1}\cdots\alpha_t^{b_t}-1,
\]
such that $ U $ is non-zero. Then the following inequality holds:
\[
|U|>\exp\left( C(t) d_\mathbb{F}^2(1+\log d_\mathbb{F})(1+\log B)A_1\cdots A_t\right),
\]
where the constant $ C(t) $ is given by
\[
C(t):= -1.4\cdot 30^{t+3}\cdot t^{4.5}, 
\]
and  $ B $ is defined as $ B\geq \max\{|b_1|,\ldots,|b_t|\}.  $ Additionally, each $ A_i $ satisfies
\[
A_i\geq \max\{d_\mathbb{F}h(\alpha_i),|\log \alpha_i|, 0.16\}, \quad \text{for} \quad  i=1,\ldots,t.
\] 
\end{theorem}
Let $ \varsigma $ be a real number, and define the distance from $ \varsigma $ to the nearest integer by $ ||\varsigma||=\min\{ | \varsigma -n | :  n \in \mathbb{Z} \}. $ 
\begin{lemma} \label{reduction}\cite[Lemma 4]{BL13}      
Let $\tau$ be an irrational and  $M$ be a positive integer. Suppose that $p/q$ is a convergent of  the continued fraction of $\tau$ such that $q>6M.$ Suppose also $K, L$ and $\mu$ are real numbers with $K>0$ and $L>1 $. If $\epsilon:=||\mu q||-M||\tau q|| >0$, then there is no positive integer solutions $ s, t, $ and $\gamma$  to the inequality
\[
0< \lvert s\tau-t+\mu \rvert <KL^{- \gamma},
\]
under the conditions
\[
s \leq M \quad\text{and}\quad \gamma \geq \frac{\log(Kq/\epsilon)}{\log L}.
\]
\end{lemma}

\begin{lemma} \label{Bang} \cite{Bang} 
Let $ n \geq 2 $ be an integer. If  $ n \neq 6 $ then $ 2^n-1 $ has a prime divisor not dividing $ 2^k-1 $ for any $ k<n. $
\end{lemma}

\section{Some Properties of $k-$Fibonacci Numbers}

The characteristic polynomial of $ k- $generalized Fibonacci sequence is given by
$$\psi_k(t)=t^k- \sum_{j=1}^k t^{k-j}.$$ 
$ \psi_k(t) $ is an irreducible polynomial over $ \mathbb{Q}[t] $ and it has exactly one real root $ \alpha(k) $ located in the outside of the unit circle \cite{22,23,26}.  The other $ k-1 $ roots are properly within the unit circle \cite{23}. The root $ \alpha(k) ,$ denoted simply as $ \alpha, $ lies in the interval $ 2(1-2^{-k})< \alpha < 2 $ and furthermore
\begin{equation}
\label{F1}
F_n^{(k)} \in  \left[ \alpha^{n-2} , \alpha^{n-1} \right]
\end{equation}
for all $n>0$ and $ k \geq 2 $ \cite{15}. Let us define
\begin{equation*}\label{fk}
f_k (t) =  \dfrac{t - 1}{2 + (k+1)(t - 2)}.
\end{equation*}
$ f_k (t)  $ has the following properties
\begin{equation}\label{fkprop}
0.5<f_k (\alpha) < 0.75  \quad \text{and} \quad \left \lvert f_k( \alpha_{(i)} ) \right\rvert <1, \quad 2 \leq i \leq k,
\end{equation}
where $ \alpha:=\alpha_{(1)}, \alpha_{(2)}  \cdots , \alpha_{(k)}  $ are all the roots of $ \psi_k(t) $  \cite[Lemma 2]{13}. Moreover, from the first one of the above inequalities, we see that $  f_k(\alpha) $ is not an algebraic integer. Moreover, in the same lemma, $ \forall k \geq 2  $, the logarithmic height of $ f_k (\alpha) $ is also calculated as
\begin{equation}\label{hfk} 
h( f_k (\alpha) ) < 3 \log { k } .
\end{equation} 
Another important property of $ F_n^{(k)}  $  for our purpose is given by
\begin{equation}\label{kF}
F_n^{(k)} = f_k (\alpha) \alpha^{n-1} + \sum_{i=2}^k  f_k ( \alpha_{(i)} ) (\alpha_{(i)})^{n-1} \quad \text{and} \quad \left \lvert F_n^{(k)} - f_k (\alpha) \alpha^{n-1} \right \rvert < 1/2,
\end{equation}
for all $ k\geq 2 $ \cite{17}. Thus,
\begin{equation}
\label{kFe}
F_n^{(k)} =  f_k (\alpha)  \alpha^{n-1} + \xi_k(n)
\end{equation}
where $ \lvert \xi_k(n) \rvert <1/2. $
\section{Proof of Theorem \ref{main1}}

Assume that \eqref{Fnml} holds. We can find the number of digits of $ F_{l}^{(k)} $ as $ d=\lfloor \log_{10}F_{l}^{(k)} \rfloor +1$  where $ \lfloor \xi\rfloor $ is the greatest integer less than or equal to $ \xi $, the floor function.  So
\begin{align*}  
d=\lfloor \log_{10}F_{l}^{(k)} \rfloor +1 \leq  1+\log_{10}F_{l}^{(k)} & \leq 1+ \log_{10}{ \alpha^{l-1}} < 1+ {(l-1)}\log_{10}{ 2 }  \\
 & < 1+ \dfrac{l-1}{3} = \dfrac{l+2}{3}
\end{align*}
and
$$  d=\lfloor \log_{10}F_{l}^{(k)} \rfloor +1 > \log_{10}F_{l}^{(k)}  \geq  \log_{10}{ \alpha^{l-2} } \geq {(l-2)} \log_{10}{ \alpha }> \frac{l-2}{5} .$$
where we used the fact that $ \log_{10}{ 2 } \cong 0.301 < 1/3 $ and   $ (1/5)<   0.208... \cong \log_{10}{ \dfrac{1+\sqrt{5}}{2} } \leq \log_{10}{ \alpha } .$ Hence we have that
\begin{equation}
\label{d}
\frac{l-2}{5} < d < \frac{l+2}{3}.
\end{equation}
Moreover, from the definition of $ d, $ we write
$$ F_{l}^{(k)} <10^d  \leq 10 F_{l}^{(k)}.$$
So, from \eqref{F1} and \eqref{Fnml}, we write
$$\alpha^{n-2} \leq F_{n}^{(k)}= 10^d F_{m}^{(k)}  +F_{l}^{(k)} \leq 10 F_{m}^{(k)} F_{l}^{(k)}+F_{l}^{(k)} \leq 11 F_{m}^{(k)} F_{l}^{(k)} \leq 11 \alpha^{m+l-2} < \alpha^{m+l+4} $$
and
$$ \alpha^{n-1} \geq F_{n}^{(k)} = 10^d F_{m}^{(k)}  + F_{l}^{(k)} > F_{m}^{(k)} F_{l}^{(k)}+F_{l}^{(k)} \geq  F_{m}^{(k)} F_{l}^{(k)} \geq \alpha^{m+l-4}.$$
Hence we have that
\begin{equation}
\label{n}
m+l-3 < n < m+l+6.
\end{equation}

\subsection{The Case $ n \leq k+1 $ } \label{nk1}
Suppose $ n \leq k+1. $ Then $ F_n^{(k)} = 2^{n-2}$, and hence  \eqref{Fnml} turns into
$$ 2^{n-2} = 2^{m-2} 10^{d}+2^{l-2}.$$
If $ l \leq m $ then we get that $ 2^{n-l} = 2^{m-l} 10^{d}+1 ,$ which means $ 0 \equiv 1 \mod 2, $ a contradiction. Suppose that $ m \leq l. $ This time we get that $ 2^{n-m} = 10^{d}+2^{l-m} .$ We write $ 2^{l-m} ( 2^{n-l}-1 )  = 10^{d}=2^d 5^d.$ From the Fundamental Theorem of Arithmetic, we get $ 2^{l-m}=2^{d} $ and $ 2^{n-l} - 1 = 5^{d} .$ From Lemma \ref{Bang}, we see that there is no positive integer solution of the last equation. Hence, there is   no solution of \eqref{Fnml} when  $ n \leq k+1. $  So, from now on, we take $ n \geq k+2. $

\subsection{An inequality for $ n $ in terms of $ k $} 
 
Now assume that $ n \geq k+2 .$ Before the detailed calculations, by a short computer program, we search all the variables satisfying \eqref{Fnml} in the range $ 0 \leq m, l <200 $ and we find only $ k-$Fibonacci numbers given in Theorem \ref{main1}. So, we assume that $ \max\{m, l\}>200. $

By using \eqref{kFe}, we rewrite \eqref{Fnml} as
$$ f_k(\alpha) \alpha^{n-1}+\xi_k(n)=( f_k(\alpha) \alpha^{m-1}+\xi_k(m))10^d+ F_l^{(k)}  ,$$
and, we get
$$ f_k(\alpha) \alpha^{n-1} - f_k(\alpha) \alpha^{m-1} 10^d = \xi_k(m)10^d+ F_l^{(k)} - \xi_k(n)  ,$$
\begin{align*}
\left\lvert  1- \dfrac{10^d}{\alpha^{n-m}} \right \rvert  & = \dfrac{1}{f_k(\alpha) \alpha^{n-1}} \left\lvert {\xi_k(m)10^d} + {F_l^{(k)}}  + {\xi_k(n)}  \right \rvert \\  
& \leq \dfrac{10 F_l^{(k)} }{\alpha^{n-1}} + \dfrac{2 \alpha^{l-1}}{ \alpha^{n-1}} + \dfrac{1}{\alpha^{n-1}} < \dfrac{14}{\alpha^{ n-l }}
\end{align*}
Hence we write
\begin{equation}
\label{lam1}
U_1 := \left\lvert  1- \dfrac{10^d}{\alpha^{n-m}}   \right \rvert <  \dfrac{14}{\alpha^{ n-l }}.
\end{equation}
Let $ \eta_1:=\alpha,$ $\eta_2:=10$  and $ b_1:=-(n-m),$  $b_2:= d .$  $ \eta_1$ and $ \eta_2 $ are belongs to the field $ \mathbb{F}=\mathbb{Q}(\alpha) $ with degree $ d_\mathbb{F}=k .$ We take $ B:=n-m \geq d. $ Indeed, in the case $ n-m < d $  from \eqref{n}, we get $l-3 <n-m < d  < (l+2)/3,$ that is $ l<4. $ Therefore we find $ d=1 $ and hence $ n=m, $ a contradiction. Since  $h(\eta_1)=({1}/{k})\log { \alpha }, $ $  h(\eta_2)= \log 10$  we take $ A_1 =  \log \alpha $ and $ A_2 = k \log { 10}.$ We have also $ U_1 \neq 0.$ If $ U_1=0 $ then we get that
$  \alpha^{n-1}=10^d $ which is a contradiction since, because of its minimal polynomial, $ \alpha $ and hence $\alpha^{n-1}  $ is a unit in the ring of integers of $\mathbb{Q}(\alpha) $ but clearly $ 10^d  $ is not.
So, by Theorem \ref{Matveev}, we get that
$$
\log{ \lvert  U_1 \rvert } > -1.4 \cdot 30^5 \cdot 2^{4.5} \cdot k^3 \log \alpha \cdot \log 10 (1+\log k)(1+\log {(n-l)}) .
$$ 
On the other hand, from \eqref{lam1}, we may directly write $ \log{ \lvert  U_1 \rvert } <  \log 14 - (n-m) \log \alpha. $ By using the fact that for all $ k \geq 3 $ and $ n \geq 5 ,$ $ 2 \log { k } > 1+\log { k }  $ and $ 2 \log {(n-m)}> (1+\log {(n-m)}),$  we obtain
\begin{equation}
\label{nlfirst}
n-l < 7.1 \cdot 10^{9} \cdot k^3 \log{k} \log {(n-m)}.
\end{equation}
Now, we revisit \eqref{Fnml} and by \eqref{kF}, we write it as
$$ f_k(\alpha) \alpha^{n-1}+\xi_k(n)= F_{m}^{(k)}10^d +  f_k(\alpha) \alpha^{l-1}+\xi_k(l).$$
Thus, we get
$$  f_k(\alpha) \alpha^{n-1} (1-  \alpha^{l-n} ) - F_{m}^{(k)}10^d  = -\xi_k(n)+ \xi_k(l).$$
\begin{align}
\label{lam2}
 U_2  := \left \lvert 1 - \dfrac{ F_{m}^{(k)}10^d  }{ \alpha^{n-1} f_k(\alpha)  (1-  \alpha^{l-n} )}  \right \rvert & < \dfrac{1}{ f_k(\alpha) (1-  \alpha^{l-n} ) } \dfrac{1}{\alpha^{n-1}} < \dfrac{6}{\alpha^{n-1}}
\end{align}
where we used the facts $ \dfrac{1}{f_k(\alpha)}<2 $ and $ \dfrac{1}{  1-  \alpha^{l-n} } < 3 $ for all $ k \geq 3 $ and $ n-l \geq 1. $

Assume that $ U_2 = 0.$ Then we have that
$$ F_{m}^{(k)}10^d = f_k(\alpha) \alpha^{n-1} (1-  \alpha^{l-n} ) .$$
By considering the image of both sides under the any one of the non-trivial isomorphisms  $ \sigma_i : \alpha \rightarrow \alpha_i $ for any $ i \geq 2 $ and by taking into account the absolute values of them, from \eqref{fkprop}, we find that
$$ 10 \leq \lvert F_{m}^{(k)}10^d \rvert = \lvert f_k( \alpha_{(i)} ) \rvert  \lvert \alpha_{(i)}^{n-1} \rvert  \lvert (1-  \alpha_{(i)}^{l-n} ) \rvert <2,$$
which is false. So, $ U_2 \neq 0 .$

Let $ \eta_1$ and $\eta_2$ be as before and $\eta_3:=\dfrac{F_{m}^{(k)}}{ f_k(\alpha) (1-  \alpha^{l-n} )}    $     with $ b_1:= -(n-1)$, $b_2:=d,$  $b_3:=1. $  All $ \eta_i$ are belongs to the number field $ \mathbb{F}=\mathbb{Q}(\alpha) \subset \mathbb{R}.$  From \eqref{hfk}, we find
$$
h(\eta_3) \leq  h( F_{m}^{(k)} ) +  h(f_k(\alpha)) + h( 1-  \alpha^{l-n} ) \leq  (m-1) \log \alpha +  3 \log k +   \dfrac{ \lvert n-l \rvert }{ k }\log \alpha +\log 2 . $$
Since $ m-1 < n-l+2 ,$  we write
$$ h(\eta_3) \leq 2 ( 2 \log k +  (n-l+1) \log \alpha  ) .$$
So $A_3= 2k (2 \log k   + (n-l+1) \log \alpha   ).$ We apply Theorem \ref{Matveev} to \eqref{lam2} and get that
\begin{equation} \label{nl2}
n-1 < 2.7 \cdot 10^{12} k^4 \log k \log(n-1) ( 2 \log k  + (n-l+1) \log \alpha  ).
\end{equation}

Now, we examine the cases $ l \leq m  $ and $ m < l $ separately. The first case implies that  $ n-m < n-l $ and hence we  write \eqref{nlfirst} as
\begin{equation*} 
n-l < 7.1 \cdot 10^{9} \cdot k^3 \log{k} \log {(n-l)}.   
\end{equation*}
It is time to cite \cite[Lemma 7]{Guzman}.
\begin{lemma}\label{Guzman} 
If $ e \geq 1 $ and $ H>(4e^2)^e $ then
$$
\dfrac{f}{(\log(f))^e} < H \Rightarrow f < 2^e H(\log(H))^e.
$$
\end{lemma}
We take $ e=1 $ and $ H:= 7.1\cdot 10^{9} k^3 \log k ,$ so that
\begin{align*}
\log H & < \log(7.1) +9\log { 10}+3\log { k }+\log( \log { k }) \\
& < 2 \log k + 27 \log k+3\log { k }+\log( \log { k })  < 35 \log k.
\end{align*}
So, by Lemma \ref{Guzman}, we obtain
\begin{equation*}\label{nlfirst2}
n-l< 4.98 \cdot 10^{11} k^3 \log^2 k.
\end{equation*} 
Since, from \eqref{n}, $ n<m+l+6 <2m+6 $ and $ m<n-l+3 $ we find that
\begin{equation*}
n< 10^{12} k^3 \log^2 k.
\end{equation*} 
Assume that $ m<l. $ Since $ \alpha < k, $  we may take $ 2 \log k + (n-l+1) \log \alpha  < (n-l+3) \log k $ and therefore by substituting the value of $ n-l $ given in \eqref{nlfirst} into \eqref{nl2}, we get
$$ n-1 < 2\cdot 10^{22} k^7 \log^2 (k) \log^2(n-1). $$
We apply  Lemma \ref{Guzman} one more time by taking $ e=2 $ and $ H= 2\cdot 10^{22} k^7 \log^2 (k) $ we find $ \log H < 80 \log k $ and hence we arrive
\begin{equation}\label{nkk}
n <  5.2 \cdot10^{26} k^7 \log^4 k.
\end{equation} 
Thus, whether $ m<l $ or not, this bound is valid. Firstly we consider a relatively small values of $ k $  that is the case $ k < 420 ,$ and then later $ k \geq 420 $.

\subsection{The Case $ k \leq 420 $}

Assume that $ 3 \leq k \leq 420. $  Let
$$ V_1 := d \log 10 - (n-m) \log \alpha .$$
Then 
$$ U_1 := \left \lvert \exp(V_1) -1 \right \rvert < 14/\alpha^{n-l}.$$
We claim that $n-l < 168.$ Suppose that $ n-l >10.$  Then  $ 14/\alpha^{n-l} < 1/2 $ and therefore $\lvert V_1 \rvert <28 / \alpha^{n-l} .$  So
\begin{equation*}
\label{g1}
0 < \left \lvert \dfrac{\log \alpha }{\log 10 } - \dfrac{d}{n-m}  \right \rvert < \dfrac{28}{  \alpha^{ n-l} (n-m) \log 10}.
\end{equation*}
If $ \dfrac{28}{ (n-m) \alpha^{ n-l} \log 10 } < \dfrac{1}{2 (n-m)^2} $ then we infer that $ d /{n-m} $ is a convergent of continued fraction expansion  of irrational ${\log \alpha}/{\log10}  ,$ say $\dfrac{p_i}{q_i}.$ Since  $ q_i \leq n-m \leq n-1 < 5.2 \cdot10^{26} k^7 \log^4 k,$ for each $ k $ we find an upper bound of indices $ i=i_0 $ and therefore the maximum value of  $a_{max}:=\max [a_0,a_1,a_2,a_3,a_4,\ldots, a_{i_0}]$ of the continued fraction of  $  \log \alpha /\log 10. $ By using a property of continued fractions, see for example \cite[Theorem 1.1.(iv)]{Hen},  we write
$$  \dfrac{1}{(a_{max}+2) \cdot (n-m)^2 } \leq \dfrac{1}{(a_i+2)(n-m)^2} < \left \lvert \dfrac{\log \alpha}{\log10} - \dfrac{d}{n-m} \right \rvert <   \dfrac{28}{\alpha^{ n-l} (n-m)\log10}  .$$
Thus, we get the inequality
$$ \alpha^{ n-l} < \dfrac{28 \cdot (a_{max}+2) \cdot 5.2 \cdot10^{26} k^7 (\log k)^4}{\log 10}, $$
which gives an upper bound for $ n-l, $ and none of them greater than 168.

If $ \dfrac{28}{ (n-m) \alpha^{ n-l} \log 10 } \geq \dfrac{1}{2 (n-m)^2} $ then we find a more strict bound for $ n-l $ as
$$ \alpha^{ n-l} < \dfrac{28 \cdot 5.2 \cdot 10^{26} k^7 (\log k)^4}{\log 10}.$$
We apply this procedure for each $ k \in [3, \cdots , 420 ] $ and we find that $ n-l <170. $ Since $ m < n-l+3 ,$ we also find a bound for $m  $ as $ m<175. $ Since $ \max\{ m, l \} >200 ,$ we may assume $ m <l .$  Let
$$ V_2 := d \log 10  - (n-1) \log \alpha + \log \left( \dfrac{F_m^{(k)}}{f_k(\alpha) (1-\alpha^{l-n})}  \right) ,$$
so that 
$$ U_2 := \left \lvert \exp(V_2) -1 \right \rvert < 6/\alpha^{n-1} <1/2.$$
Hence, we get that
\begin{equation}
\label{g2}
0 < \left \lvert d \dfrac{\log 10 }{\log \alpha }  - (n-1) + \dfrac{  \log \left( \dfrac{F_m^{(k)}}{f_k(\alpha) (1-\alpha^{l-n})}  \right)  }{\log 10 }    \right \rvert < \dfrac{12}{\alpha^{n-1} \log 10  }.
\end{equation}
We calculate the value of $\epsilon_{(k, m, n-l)} := ||\mu_{(k, m, n-l)} q_{i} || -M_k || \tau_k q_{i} ||  $ for each $ m, n-l \in \{1, 2, \ldots, 175 \} ,$   where $ \tau_k = \dfrac{\log 10 }{\log \alpha } $ and
$$ 
\mu_{(k, m, n-l)}  := \dfrac{  \log( \dfrac{F_m^{(k)}}{f_k(\alpha) (1-\alpha^{l-n})}  )  }{\log 10 } .
$$
In the case of $\epsilon_{(k, m, n-l)}< 0 $ for any values of $m  $ or $ n-l ,$ we repeat the same procedure by using $ q_{i+1} $ instead of $ q_i ,$ and for each $ k ,$ we obtain an appropriate $\tau_k$ such that $ \epsilon_{(k, m, n-l)} > 0.$ Thus, the conditions on Lemma \eqref{reduction} are satisfied and hence we apply  Lemma \ref{reduction} to the equation \eqref{g2}. As a result, for each $ 3 \leq k \leq 420,$  we obtain an upper bound for $ n-1 $  and none of them exceed 171. Some of these bounds are 135 for $ k=3 $ and $ k=10, $ 164 for $ k=100 ,$  $ k=200, $ $ k=300 $ and $ 167 $ for $ k=400. $ Taking into account these bounds of $ n-1 $ together with \eqref{d}, we write a short computer programme to check all variables in the range $ 3 \leq k \leq 420, $  $ 0 \leq m < n-l+3 ,$  $ k+2 \leq n \leq n(k)$ satisfying  \eqref{Fnml}. We find only the results stated in Theorem \ref{main1}. Now we need to examine the case $ k>420. $

\subsection{The Case $ k > 420 $}
\begin{lemma} \label{Guzman2} \cite[Lemma 3]{14}
Let $ n<2^{k/2} .$ Then $F_n^{(k)}  $ satisfies the estimate
$$  F_n^{(k)}=2^{n-2}(1+\varsigma(n,k)), \quad \text{where} \quad \lvert \varsigma(n,k) \rvert  < \dfrac{2}{2^{k/2}} .$$
\end{lemma}
For $ k>420 ,$ the inequality 
$$ n < 5.2 \cdot 10^{26} k^7 (\log { k })^4 < 2^{k/2}, $$ 
holds and hence from the above Lemma \ref{Guzman2}, we may write
\begin{equation*}
F_n^{(k)}=2^{n-2} + 2^{n-2} \varsigma(n,k) \quad \text{and} \quad F_m^{(k)}=2^{m-2} + 2^{m-2} \varsigma(m,k).
\end{equation*}
Thus, substituting these values of $ F_n^{(k)}$ and $ F_m^{(k)} $ into \eqref{Fnml}, we get 
\begin{align*}
\left \lvert  2^{n-2} - 10^d 2^{m-2} \right \rvert = & \left \lvert  10^d 2^{m-2} \varsigma(m,k) + F_l^{(k)} - \varsigma(n,k)2^{n-2} \right \rvert \\
& < \dfrac{2^{m-1} 10 \alpha^{l-1} }{2^{k/2}} + \alpha^{l-1} + \dfrac{2^{n-1} }{2^{k/2}}.
\end{align*}
Thus, an easy calculation yields
\begin{align*}
\left \lvert 1- \dfrac{10^d}{2^{ n-m }}  \right \rvert  & \leq \dfrac{10}{ 2^{ n-l-m } 2^{k/2}} +\dfrac{1}{2^{ n-l-1 }} + \dfrac{2}{2^{k/2}} \\
& \leq \dfrac{10 \cdot 2^3+2}{2^{k/2}} + \dfrac{1}{ 2^{ n-l-1 } } < \dfrac{1}{2^{{k/2}-7}} + \dfrac{1}{2^{ n-l-1 }} 
\end{align*}
since, from \eqref{n}, $ \dfrac{1}{2^{ n-l-m }} < 2^3.$ Therefore, we write
\begin{equation*} \label{lam3}
 U_3 := \left \lvert 1- \dfrac{10^d}{2^{ n-m }} \right \rvert < \dfrac{1}{2^\lambda} ,
\end{equation*}
where $ \lambda :=\min\{ (k/2)-8, n-l-2 \} .$ By taking  $ \eta_1:=10,$ $\eta_2:=2,$  and $ b_1:=d,$  $b_2:=-(n-m),$ we get $ \mathbb{F}=\mathbb{Q} ,$ $ d_\mathbb{F}=1 ,$ $ B: = n-m.$ Clearly $ U_3 \neq 0. $  So, from Theorem \ref{Matveev}, we get that
$$
\lambda \log 2 < 1.4 \cdot 30^5 \cdot 2^{4.5}  (1+\log {(n-m)}) \log {10} \cdot \log 2 .
$$
So, we find that
\begin{equation*}
 \lambda  <  3.6 \cdot 10^{9} \cdot \log {(n-m)}.
\end{equation*}
Note that
\begin{align*}
\log {(n-m)}  & \leq  \log {(n-1)} <  \log ( 5.2 \cdot 10^{26} k^7 \log^4 { k } ) \\
& < 14 \log(420)+7 \log k +4 \log \log { k } < 26 \log k.
\end{align*}
Therefore, we have that
\begin{equation*}
 \lambda  <  9.37 \cdot 10^{10} \cdot \log {k}.
\end{equation*}
A quick computation shows that $ \lambda:=(k/2)-8 $ leads to $ k < 10^{15}.$
If $ \lambda: = n-l-2  $ then we get 
\begin{equation*} \label{d37}
n-l < 9.4 \cdot 10^{10} \log k. 
\end{equation*}
We turn back to \eqref{Fnml} and using the estimates given in Lemma \ref{Guzman2} we rewrite it as
\begin{equation*} 
2^{n-2} - 2^{m-2} 10^d - 2^{l-2} =    2^{m-2} \varsigma(m, k) 10^d + 2^{l-2} \varsigma(l, k)  - 2^{n-2} \varsigma(n, k).
\end{equation*}
\begin{equation*} 
2^{n-2} (1 - 2^{l-n} ) - 2^{m-2} 10^d  =    2^{m-2} \varsigma(m, k) 10^d + 2^{l-2} \varsigma(l, k) - 2^{n-2} \varsigma(n, k).
\end{equation*}
Hence we obtain
\begin{equation*}\label{lam4}
 U_4  =\left \lvert  1- \dfrac{10^d}{2^{n-m} (1 - 2^{l-n} ) } \right \rvert < \dfrac{168}{2^{k/2}},
\end{equation*}
where we used the following two facts
\begin{equation*}
\dfrac{1}{1 - 2^{l-n}} < \dfrac{1}{1-1/2} =2,
\end{equation*}
\begin{equation*}
\dfrac{2^{m-2} 10 \cdot 2^{l-1}}{2^{n-2}} +\dfrac{2^{l-2}}{2^{n-2}} +1 < \dfrac{5}{2^{n-m-l}}+1+1<\dfrac{5}{2^{-3}}+2 = 42
\end{equation*}
Let  $ \eta_1:=10,$ $\eta_2:=2,$ $ \eta_3:= (1 - 2^{l-n} )^{-1}  $ and $ b_1:d,$  $b_2:=-(n-m) ,$ $b_3:=1. $ Then $ \mathbb{F}=\mathbb{Q} ,$ $ d_\mathbb{F}=1 ,$ $ B: = n-m < 26 \log k.$  $ A_1 = \log 10, $ $ A_2 = \log 2 $ and
$$  h(\eta_3) \leq \lvert n-l \lvert \log 2+ \log 2 =(n-l+1) \log 2 < 6.6 \cdot 10^{10} \log k $$

Exactly the same argument as in the subsection \ref{nk1} shows that $ U_4 \neq 0 .$  Thus, application of Theorem \ref{Matveev} to $ U_4 $ and the fact that  $ \log {U_4}< \log{168}-(k/2) \log 2 $ gives $ k <1.2 \cdot 10^{24} \log^2 {k}.$ So, by Lemma \ref{Guzman}, we get
\begin{equation*}
k<10^{28}.
\end{equation*}
Thus, by \eqref{nkk}, this bound of $ k $ leads to  bound of $ n $ as
\begin{equation}\label{n229}
n < 9 \cdot 10^{229}.
\end{equation}
\subsection{Reducing The Bound on k}
Let 
\begin{equation}\label{g3}
V_3 := d \log 10 - (n-m) \log 2
\end{equation}
Then $ U_3 :=\left \lvert \exp(V_3)-1 \right \rvert < \dfrac{1}{2^\lambda} .$
Suppose that $ \lambda>2 $. Then, $ \dfrac{1}{2^\lambda} < \dfrac{1}{2}$ and hence we get that $ \lvert V_3 \rvert < \dfrac{2}{2^\lambda} .$ Then, from \eqref{g3}
$$ \left \lvert \dfrac{\log2}{\log10} - \dfrac{d}{n-m}  \right \rvert < \dfrac{2}{ (n-m)   2^\lambda  \log10}.   $$
First assume that $ \dfrac{2}{2^\lambda (n-m)\log10} < \dfrac{1}{2 (n-m)^2} .$ Then, this inequality means $ \dfrac{d}{n-m} $ is a convergent of continued fractions of ${\log2}/{\log10}  ,$ say $\dfrac{p_i}{q_i}.$ Since $ \gcd( p_i ,  q_i)=1 $, we deduce that $ q_i \leq n-m \leq n-1 < 9 \cdot 10^{229}. $ A quick calculation shows that $ i<468. $ Let  $[a_0,a_1,a_2,a_3,a_4,\ldots]=[0,3,3,9,2,  \ldots]$ be the continued fraction expansion of  $  \log 2 /\log 10. $ 
Then  $ \max\{a_i\}  = 5393 $ for $ i=0,1,2, \ldots , 468. $ So,  we write
$$  \dfrac{1}{5395 \cdot (n-m)^2 } \leq \dfrac{1}{(a_i+2)(n-m)^2} < \left \lvert \dfrac{\log2}{\log10} - \dfrac{d}{n-m} \right \rvert <   \dfrac{2}{  (n-m) 2^\lambda \log10}  .$$
Thus, we have that
$$ 2^\lambda < \dfrac{2 \cdot 5395 \cdot 9 \cdot 10^{229}}{\log 10} < 4.3 \cdot 10^{233}, $$ 
that is $ \lambda < 777.  $  On the other hand, the inequality 
$$ \dfrac{1}{2 (n-m)^2} \leq \dfrac{2}{2^\lambda (n-m)\log10}  $$
implies
$$ 2^\lambda < \dfrac{4(n-m)}{\log10} < \dfrac{4 \cdot 9 \cdot 10^{229}}{\log 10}  < 1.6 \cdot 10^{230} < 2^{765}. $$
So we see that $ \lambda < 777 $ holds in this case too.  If $ \lambda=k/2 -8$ then $ k <1570.$ If $ \lambda = n-l-2 $ then $ n-l < 780.$ Let
\begin{equation*}
V_4 = \left \lvert  d \log 10 - (n-m) \log 2  -  \log( 1-2^{l-n} ) \right \rvert .
\end{equation*} 
Then
\begin{equation*}
U_4:= \lvert \exp( V_4 )-1 \rvert < \dfrac{168}{2^{k/2}}<\dfrac{1}{2}
\end{equation*}
So
\begin{equation*}
0< \left \lvert \dfrac{V_4 }{|\log 2}  \right \rvert < \left \lvert d \dfrac{\log 10 }{\log 2 } - (n-m) -  \dfrac{ 1-2^{l-n}    }{{\log 2}} \right \rvert  < \dfrac{336}{2^{k/2}{\log 2 }}
\end{equation*}
We take $ M:= 9 \cdot 10^{229} >n-1 \geq n-m >d $ and  $\tau = \dfrac{ \log{10}}{\log 2}. $  Note that $q_{472} > 6M.$ So, we take

$\epsilon_{n-l}:= ||\mu_{n-l} q_{472} || -M || \tau q_{472} ||  $ for each $ n-l \in \{3, 4, \ldots, 780 \} $ where
$$ \mu_{n-l} :=  - \dfrac{ \log{ ( 1-2^{l-n} )   }}{{\log 2}}  $$
and we find that $0.000957 \leq \epsilon_{n-l}  $ for  all $n-l .$ Let $ K:= \dfrac{336}{\log 2} ,$ $ L:=2 $ and $ \omega :=k/2 .$ By a simple calculation, we see that
$$ \dfrac{ \log{ \left( Kq_{472}/ 0.000957 \right) }  }{\log L} < 795,$$
and hence from Lemma \ref{reduction}  we  find $ k < 1590 .$ Hence, from \eqref{nkk},  $ n < 4 \cdot 10^{52} .$ We repeat the same procedure starting with this subsection and  find that $ \lambda < 188. $ Thus $ \lambda = k/2-8 $ leads to $ k <385, $ which is a contradiction by choosing the range of $ k$ whereas $ \lambda = n-l-2 $ gives $ n-l < 190.$ So, using this new upper bound $ n-l < 190,$  we work on $ V_4 $ as we did before but with $ M:=4 \cdot 10^{52} $ and $  q_i :=q_{120} .$ We find that
$$ 0.001034 \leq  \epsilon_{n-l}  $$ for  all $n-l .$ 
By Lemma \ref{reduction}, we find $ k/2 < 205,$ that is $ k < 410 $
which is the desired contradiction to our assumption that $k > 420.$ Thus we conclude that \eqref{Fnml} has no positive integer solution when $ k >420. $ This concludes the proof.


\begin{thebibliography}{99}                                                                                               %


\bibitem{Adedji1}
Ad{\'e}dji, K. N., Faye, M. N.,  Togb{\'e}, A. (2024). On the Diophantine equations $P_n= b^d Q_m+ Q_k$ and $Q_n= b^d P_m+ P_k$ involving Pell and Pell–Lucas numbers. Proceedings Mathematical Sciences, 134(1), 14.


\bibitem{Adedji2}
Ad{\'e}dji, K. N., Kandhil, N., Togb{\'e}, A. (2024). On $b-$concatenations of Padovan and Perrin numbers. Bolet{\'\i}n de la Sociedad Matem{\'a}tica Mexicana, 30(2), 71.


\bibitem{Adedji3}
Ad{\'e}dji, K. N., Trebje{\v{s}}anin, M. B. (2024). On Mixed $ B $-Concatenations of Pell and Pell–Lucas Numbers which are Pell Numbers. Mathematica Pannonica, 30(1), 91-104.



\bibitem{Alan}
Alan, M.: \textit{On Concatenations of Fibonacci and Lucas Numbers}, Bull. Iran. Math. Soc. \textbf{48} (2022) 2725--2741.


\bibitem{Altassan}
Altassan, A., Alan, M. (2024). Fibonacci numbers as mixed concatenations of Fibonacci and Lucas numbers. Mathematica Slovaca, 74(3), 563-576.




\bibitem{Bang}
Bang, A. S. (1886). Taltheoretiske Undersøgelser.(Fortsat, se S. 80). Tidsskrift for mathematik, 4, 130-137.

\bibitem{Banks}
Banks, W. D. Luca, F.: \textit{Concatenations with binary reccurent sequences},  Journal of Integer Sequences, \textbf{8}, 05.1.3 (2005).


\bibitem{13}
Bravo, J.J., G´omez, C.A.G., Luca, F.: Powers of two as sums of two k-Fibonacci
numbers. Miskolc Math. Notes 17, 85–100 (2016)


\bibitem{14}
Bravo, J.J., G´omez, C.A., Luca, F.: A Diophantine equation in k-Fibonacci
numbers and repdigits. Colloq. Math. 152, 299–315 (2018)




\bibitem{15}
Bravo, J.J., Luca, F.: Powers of two in generalized Fibonacci sequences. Rev.
Colombiana Mat. 46, 67–79 (2012)






\bibitem{BL13}
Bravo J.J., Luca F., {\em On a conjecture about repdigits in k-generalized Fibonacci sequences}, Publ. Math. Debrecen, {\bf 82} (2013), 623-639.




\bibitem{Bravo23}
Bravo, E. F. (2023). On concatenations of Padovan and Perrin numbers, Mathematical Communications, 28(1), 105-119.







\bibitem{17}
Dresden, G., Du, Z.: A simplified Binet formula for k-generalized Fibonacci
numbers, J. Integer Sequences 17, Article 14.4.7 (2014),


\bibitem{DP} 
Dujella A and Peth\H o A,  A generalization of a theorem of Baker and Davenport, {\em Quart. J. Math. Oxford Ser.} {\bf 49}(195) (1998) 291-306 


\bibitem{Guzman} 
S. Guzman and F. Luca, Linear combinations of factorials and S-units
in a Binary Recurrence Sequence, Ann. Math. Qu´ebec 38 (2014), 169–
188.


\bibitem{Hen}
D. Hensley, \emph{ Continued fractions}, World Scientific Publishing Co. Pte. Ltd., Hackensack, NJ, (2006)



\bibitem{Erduvan}
Erduvan, F. (2023). On Concatenations of Two Padovan and Perrin Numbers. Bulletin of the Iranian Mathematical Society, 49(5), 62.


\bibitem{Guneymerve}
G{\"u}ney Duman, M. (2024). Padovan numbers that are concatenations of a Padovan number and a Perrin number. Periodica Mathematica Hungarica, 1-16.





\bibitem{Matveev}
Matveev E M, An explicit lower bound for a homogeneous rational linear form in the logarithms of algebraic numbers, II, {\em  Izv. Ross. Akad. Nauk Ser. Mat.} {\bf 64}(2000) 125-180, {\em Translation in Izv. Math.} {\bf 64}(2000) 1217-1269


\bibitem{22}
Miles, E.P., Jr.: Generalized Fibonacci numbers and associated matrices. Am.
Math. Mon. 67, 745–752 (1960)

\bibitem{23}
Miller, M.D.: Mathematical notes: on generalized Fibonacci numbers. Am. Math.
Mon. 78, 1108–1109 (1971)





\bibitem{26}
Wolfram, D.A.: Solving generalized Fibonacci recurrences. Fibonacci Quart.
36(2), 129–145 (1998)

\end{thebibliography}
\end{document}